\documentclass[12pt]{article}
\usepackage{amsthm,amsfonts}
\newtheorem{thm}{Theorem}[section]

\newtheorem{lem}[thm]{Lemma}
\newtheorem{prop}[thm]{Proposition}

\theoremstyle{definition}
\newtheorem{defn}[thm]{Definition}
\theoremstyle{remark}
\newtheorem{rem}[thm]{Remark}
\newcommand{\ex}{\noindent {\bf Example }}
\renewcommand{\proof}{\noindent {\bf Proof: }}

\title{Test configurations and Geodesic rays}
\author{Xiuxiong Chen\footnote{Partially supported by NSF grant}, Yudong Tang}
\begin{document}
\maketitle
\tableofcontents

\section{Introduction}
The purpose of this paper is to explore the connection between geodesic rays in the space of K\"ahler metrics in algebraic manifold and test configurations~\cite{DonaldsonToric}.   This is a continuation of \cite{Chenbounded} in some  aspects.  In \cite{ChenCalabi}, the first named author and E. Calabi proved that the space of K\"ahler potentials is a non-positive curved space in the sense of Alexanderov.  As a consequence, they proved that for any given geodesic ray and any given K\"ahler potential outside of the given ray,
there always exists a geodesic ray in the sense of metric distance ($L^2$ in the K\"ahler potentials) which initiates from the given K\"ahler potential and parallel to the initial geodesic ray.  The initial geodesic ray, plays the role of prescribing an asymptotic direction for the new geodesic ray out of any other K\"ahler potential.  When the initial geodesic ray is smooth and is tamed by a bounded ambient geometry, 
 the first named author  \cite{Chenbounded} proved the existence of relative $C^{1,1}$ geodesic ray from any initial K\"ahler potential.  (These definitions can be found in Section 2.)   Similarly, as remarked in  \cite{Chenbounded},  
a test configuration should plays a similar role.  One would like to know if it induces a relative $C^{1,1}$ geodesic ray from any other K\"ahler potential in the direction of test configuration.   In \cite{Tianarezzo}, Arezzo and Tian proved a surprising result that  for a smooth test configuration with
analytic (smooth) central fiber, there always exists  an asymptotic smooth geodesic ray from
fibre which is close enough to the central fiber.  A natural question,  motivated by Arezzo-Tian's work, is if there  exists a relative geodesic ray from arbitrary initial K\"ahler metric which
also reflects the same geometry (i.e., degenerations) of the underlying test configuration.  In section 3, we prove

\begin{thm}\label{main}Every smooth test configuration induces a relative $C^{1,1}$ geodesic ray.\footnote{Following ideas of \cite{Chenbounded}, the smooth assumption can be reduced to a lower bound of the Riemannian  curvature of the total space.}
\end{thm}



Test configurations can be viewed as algebraic rays, which are geodesics in a finite dimensional subspace( with new metric) of space of K\"ahler metrics. The geodesic rays induced by a test configuration are the rays parallel to the algebraic ray. They automatically have bounded ambient geometry introduced by the first named author~\cite{Chenbounded}. 


\begin{thm} For simple test configuration~\footnote{Definition \ref{SimpleTest} }, if the induced geodesic ray is smooth regular\footnote{Definition \ref{SmoothReg}, it is also equivalent to Definition \ref{SuperReg} in this case}, then the generalized Futaki invariant agrees with the $\yen$ invariant\footnote{The $\yen$ invariant is defined by the first named author~\cite{Chenbounded}}.
\end{thm}

The  Futaki invariant was  initially introduced by Futaki~\cite{Futaki} as obstruction to
the existence of K\"ahler Einstein metrics.  E. Calabi~\cite{CalabiFutaki} generalized it to
be an obstruction for the existence of  constant scalar curvature (cscK) metrics.  It was then  generalized by Ding, Tian~\cite{DTFutaki} in the case of special degeneration.   When Tian studied the existence of K\"ahler Einstein metrics with positive scalar curvature, he  \cite{inventionpaper97}   introduced the notion
of K stability by using this generalized Futaki invariant in special degeneration.  In the same paper,
G. Tian proved that the existence of KE metric implies K semi-stability.  In 2002, S. K. Donaldson 
formulated an algebraic Futaki invariant and defined an equivalent version of  K stability on more general test configuration by using the algebraic Futaki invariant.  One important step in Donaldson's approach is to prove a theorem similar to our  Theorem 1.2 for the
generalized Futaki invariant of Ding-Tian and the algebraic Futaki invariant of Donaldson. 
\\
 
On the other hand, the generalized Futaki invariant or algebraic Futaki invariant is an algebraic notion which relates to the stability of projective manifolds. It is a well-known conjecture that the existence of constant scalar curvature metrics, or extremal K\"ahler metrics more generally, is equivalent to some kind of algebraic stability (Yau-Tian-Donaldson conjecture). In \cite{Chenbounded},  the first named author use $\yen$ invariant to define geodesic stability. Theorem 1.2 states
that geodesic stability in the algebraic manifold, is a proper generalization of K stability, at least conceptually. The first named author believes that the existence of KE metrics is equivalent to the geodesic stability introduced in \cite{Chenbounded}.  \\

The Yau-Tian-Donaldson conjecture is a central problem in K\"ahler geometry now.  Through the hard work of  many mathematicians, we now know more about one direction ( from existence to stability),
cf.  Tian\cite{inventionpaper97},  Donaldson~\cite{Donaldsonlowerbound} , Mabuchi\cite{Stability}, Paul-Tian\cite{PaulTian},  Chen-Tian\cite{ChenTian}.... But on the direction from algebraic stability to existence, few progress has been made though. However, in toric manifolds, there has been special results of Donaldson~\cite{DonaldsonToric} and Zhou-Zhu\cite{Zhu}.\\

There is an intriguing  work by V.Apostolov, D.Calderbank, P.Gauduchon, C.W.Tonnesen-Friedman \cite{Gaudachon}. They constructed an example
which is suspected to be algebraically K stable\footnote{Generalized K stable for extremal K\"ahler metrics, cf. \cite{Gabor}.}, but admits no extremal K\"ahler metric. Perhaps one
might speculate that, the geodesic stability introduced in \cite{Chenbounded} is one of the possible alternatives
since it appears to be stronger than K stability and it is a non algebraic notion in nature.\\

The converse to Theorem~\ref{main} is widely open. In other words, it is hard to compactify a geodesic ray. The  rays induced by any test configuration is very special in many aspects. For instance, the foliation of a smooth geodesic ray is not periodic in general. However, for the geodesic rays induced from a test configurations, the foliation is always periodic.  Unfortunately, having a periodic orbit does not
appear to be enough to construct a test configuration.  It would be a very intriguing problem to find
 a sufficient condition so that we can ``construct" a test configuration out of a ``good" geodesic ray. \\

\noindent {\bf Question A} {\it 
Is there a canonical method to construct some test configuration/algebraic ray such that it reflects the same degeneration of a geometric ray? What is natural geometric conditions on the ``good" geodesic ray? }\\

Our second main result is to establish the correspondence between smooth regular solutions of
Homogeneous complex Monge Ampere equation (HCMA) on simple test
configurations and some family of holomorphic discs in an ambient space $\mathcal{W}$ which will be explicitly constructed. We
prove, in section 5,

\begin{thm}There is a one to one correspondence between smooth regular
solutions of HCMA on simple test configuration $\mathcal{M}$ and
families of holomorphic discs in $\mathcal{W}$ with proper boundary
condition.
\end{thm}

Note that in the case of disc, roughly speaking,  S. K. Donaldson \cite{DonaldsonSetup} and Semmes \cite{HCMAsemmes} established first such a correspondence between
the regularity of the solution of the HCMA equation and the smoothness of the moduli space of
holomorphic discs whose boundary lies in some totally real sub-manifold.   The theorem above is a
generalization of Donaldson's result.  Following this point of view,  the regularity of the solution
is essentially the same as the smoothness of the  moduli space of these holomorphic discs
under perturbation.  As in \cite{DonaldsonSetup}, we proved the openness of smooth
regular solutions in Section 6

\begin{thm}Let $\rho(t)$ be a smooth regular geodesic ray induced by a simple test configuration. Then there exists a parallel smooth regular geodesic ray for any initial point sufficiently close to  $\rho(0)$ in $C^\infty$ sense.
\end{thm}

An immediate corollary is that the smooth geodesic ray constructed by Arezzo-Tian is open
for small deformation of the initial K\"ahler potential.  One may wonder what about the closeness of these  solutions? Note that the first named author and Tian~\cite{ChenTian}  studied  the compactness of these holomorphic discs in the disc setting and we believe that the technique of \cite{ChenTian} can be extended
over here.\\

In Section 7, as a special case, we explore the geodesic rays induced by toric degenerations~\cite{DonaldsonToric}.  In particular, we found plenty of geodesic rays whose regularity is at most  $C^{1,1}$ globally. We state a theorem with a sketch of the proof:
\begin{thm} The geodesic ray induced by a toric degeneration has the initial direction equal to the extremal function in the polytope representation.
\end{thm}

More interestingly,  we can write down the geodesic ray explicitly in polytope representation. Thus,  the  various invariants and  energies can be calculated explicitly.\\

\noindent {\bf Acknowledgments} Both authors are grateful to G. Tian for many insightful discussions. The first named author is grateful to S. K. Donaldson for many discussions in this subject.

\section{Preliminary}
\subsection{Geodesic rays in K\"ahler potential space} Let
$(M,\omega,J)$ be a compact K\"ahler manifold of complex dimension
$n$. This means $J$ is an integrable complex structure and the
symplectic form $\omega$ is compatible with $J$. In another word,
$\omega(J\cdot,J\cdot)=\omega(\cdot,\cdot)$, and
$g=\omega(\cdot,J\cdot)$ is a metric.

In local complex coordinates $z_\alpha=x_\alpha+iy_\alpha$, denote
the metric $g=\omega(\cdot,J\cdot)$ by
$g_{\alpha\bar\beta}dz^\alpha\otimes dz^{\bar\beta}$. 
$g_{\alpha\bar\beta}$ is the complexification of the real metric
$g_{ij}$.

By definition, $\omega= {\sqrt{-1}\over 2} g_{\alpha\bar\beta}dz^\alpha\wedge
dz^{\bar\beta}$. Let
\begin{equation}
\mathcal{H}=\{\phi\in
C^\infty(M):g_{\alpha\bar\beta}+\frac{\partial^2\phi}{\partial
z_\alpha
\partial z_{\bar\beta}}>0\}
\end{equation}
It follows the $\partial\bar\partial$ lemma that $\mathcal{H}$ is
the moduli space of all K\"ahler metrics in the class $[\omega]$.

$\mathcal{H}$ is an infinite dimensional manifold with formal
tangent space $T\mathcal{H}_{\phi}=C^\infty(M)$. Mabuchi defined a
metric as the following: Let $\phi_1,\phi_2\in T\mathcal{H}_\phi$.
\begin{equation}
<\phi_1,\phi_2>_{\omega_\phi}=\int_M \phi_1\phi_2
d\mu=\int_M\phi_1\phi_2\frac{\omega_\phi^n}{n!}=\int_M\phi_1\phi_2\frac{(\omega+i\partial\bar\partial\phi)^n}{n!}
\end{equation}

Under this metric, the
geodesic equation for curve $\phi(t)\in\mathcal{H}$ is the following:
\begin{equation}
\ddot\phi-g_\phi^{\alpha\bar\beta}\dot\phi_{\alpha}\dot\phi_{\bar\beta}=0
\end{equation}
It is just the Euler-Lagrange equation of the energy
$E(\phi(t))=\int_0^1\int\dot\phi^2\frac{\omega_\phi^n}{n!}dt$.
Donaldson and Semmes transformed the geodesic equation into a
Complex Monge-Ampere equation: Let $\Sigma=[0,1]\times S^1$, a
Riemann surface. Now $\phi$ is originally defined for $t\in[0,1]$. 
Extend $\phi$ to be $S^1$ invariant function on $\Sigma$. Let
$z=t+is$ be complex coordinate of $\Sigma$, $w_\alpha$ be local
coordinates on $M$. Then the geodesic equation is transformed into
\begin{equation}
\det\left(%
\begin{array}{cc}
  g_{\alpha\bar\beta}+\phi_{\alpha\bar\beta} & \phi_{\alpha\bar z} \\
  \phi_{z\bar\beta} & \phi_{z\bar z} \\
\end{array}%
\right)=0
\end{equation}
In another word, it is $(\Omega+i\partial\bar\partial\phi)^{n+1}=0$
on $M\times \Sigma$, where $\Omega=\pi^*\omega$ is the pull back of
$\omega$ by the projection $\pi:M\times \Sigma\rightarrow M$.

Now, the geodesic connecting two points $\phi_0$ and $\phi_1$ is
the solution of:
\begin{eqnarray}
\det\left(%
\begin{array}{cc}
  g_{\alpha\bar\beta}+\phi_{\alpha\bar\beta} & \phi_{\alpha\bar z} \\
  \phi_{z\bar\beta} & \phi_{z\bar z} \\
\end{array}%
\right)&=&0\texttt{    on    }M\times\Sigma\\
\phi&=&\phi_0\texttt{    on    }M\times {0}\times S^1\\
\phi&=&\phi_1\texttt{    on    }M\times {1}\times S^1
\end{eqnarray}

\begin{defn}\label{SmoothReg}Smooth Regular solution: We call $\phi$ a smooth regular solution (sometimes smooth solution for simplicity) of the Monge-Ampere equation, if $\phi$ is smooth and $g_{\alpha\bar\beta}+\phi_{\alpha\bar\beta}>0$ on fibers.
\end{defn}

In ~\cite{ChenC11}, The first named author proved the existence of a $C^{1,1}$
solution to above equation. He used the continuity method to solve $\det=\epsilon f$ equation, and proved the following: For every $\epsilon>0$, there is a unique smooth solution $\phi_\epsilon$ with $|\partial\bar\partial\phi_\epsilon|<C$. The $C$ only depends on the background metric and the manifold. In fact, his proof works for Monge Ampere equation on general
compact complex manifold with boundary. He also proved the
uniqueness of the limit when $\epsilon\rightarrow 0$. Notice that the uniqueness is expected since $\mathcal{H}$ is negatively curved space. Donaldson~\cite{DonaldsonSym} showed $\mathcal{H}$ is negatively curved in formal sense and later, the first named author and Calabi~\cite{ChenCalabi} proved it is negatively curved in the sense of Alexanderof.

The regularity beyond $C^{1,1}$ is missing. Our example in section 7 showed a solution with no global $C^3$ bound. A similar setup~\cite{DonaldsonSetup} to the geodesic equation is concerned Monege Ampere equation on $M\times D$ instead of $M\times(I\times S^1)$. In that setup, 
Donaldson showed there exists boundary value such that there is no
smooth regular solution. In this direction, a deep analytic result is~\cite{ChenTian}. The first named author and Tian characterize the singularity in detail by
analyzing the holomorphic discs associated to a solution.

In geodesic ray case, the
equation holds on $M\times [0,\infty)\times S^1$ instead of $M\times I\times S^1$. By changing
variable: $z=e^{-(t+is)}$, the strip $[0,\infty)\times S^1$ goes
to a punched disc. The equation becomes
$(\Omega+i\partial\bar\partial\phi)^{n+1}=0$ on $M\times (D-0)$.

\subsection{Test configuration and equivariant embedding}
Test configuration is defined by Donaldson~\cite{DonaldsonToric}. He used test
configurations to study the relation between stability of projective
manifolds and the existence of extremal K\"ahler metrics. Test configuration is
parallel to the notion "the special
degeneration", defined by Tian earlier. Briefly speaking, they both describe a certain degeneration of
K\"ahler manifolds. On the other hand, the geodesic ray represents the
degeneration of K\"ahler metrics. So they are naturally related.

Following Donaldson's definition,
\begin{defn} Let $L\rightarrow M$ be an ample line bundle over
a compact complex manifold. A test configuration $\mathcal{M}$
consists of:
\begin{enumerate}
  \item a scheme $\mathcal{M}$ with a $C^*-$action.
  \item a $C^*-$equivariant line bundle
  $\mathcal{L}\rightarrow\mathcal{M}$.
  \item a flat $C^*-$equivariant map $\pi:\mathcal{M}\rightarrow C$,
  where $C^*$ acts on $C$ by multiplication.
Any fiber $M_t=\pi^{-1}(t)$ for $t\neq 0$ is isomorphic to
$M$. The pair $(L^r,M)$ is isomorphic to
$(\mathcal{L}|_{M_t},M_t)$ for some $r>0$, in particular,
$(L^r,M)=(L_1,M_1)$.
\end{enumerate}
\end{defn}

Test configuration is more explicit in the view of equivariant embedding~\cite{Ross}.
Without loss of generality, assume $r=1$. For large
$k$, $\mathcal{L}^k\rightarrow\mathcal{M}\rightarrow C$ can be
embedded into $\mathcal{O}(1)\rightarrow P^N\times C\rightarrow C$
equivariantly. It means there is a $C^*$ action on
$\mathcal{O}(1)\rightarrow P^N\times C\rightarrow C$, which restricts to
the $C^*$ action of the embedded
$\mathcal{L}^k\rightarrow\mathcal{M}\rightarrow C$. In fact, the
embedding of each fiber $M_t$ is just the Kodaira embedding by the
linear system $H^0(M_t,\mathcal{L}^k|_{M_t})$. Moreover, one can
make the $S^1$ action on $\mathcal{O}(1)\rightarrow P^N\times
C\rightarrow C$ unitary.

In the rest of the paper, we always treat test configurations as
equivariantly embedded with $r=1,k=1$. Therefore, we work at a
subspace of $P^N\times C$. Also, in geodesic ray problem, there is no loss of
generality to only look at truncated test configuration
$\mathcal{M}\rightarrow D$.

At last, we define a special kind of test configuration. Geometrically speaking, it is the best behaved test configuration. 
\begin{defn}\label{SimpleTest} Simple test configuration: A test configuration $\mathcal{M}\subset
P^N\times D$ is called simple if the total space is smooth
($\mathcal{M}$ is a smooth sub-manifold of $P^N\times D$) and the
projection $\pi:\mathcal{M}\rightarrow D$ is submersion everywhere.
\end{defn}

By definition, the central fiber of a simple test configuration is automatically smooth.
\section{Relative $C^{1,1}$ geodesic ray from smooth test configuration}
\subsection{Existence}
As mentioned before, test configuration represents some degeneration of a K\"ahler manifold along a $C^*$ action. Geodesic ray represents a degeneration of K\"ahler metrics along a punched disc. So it is natural to relate the truncated test configuration to a geodesic ray.  We have the following theorem:

\begin{thm}\label{coninduceray}A smooth truncated test configuration $\mathcal{M}\rightarrow D$ induces a relative $C^{1,1}$
geodesic ray from any given initial point $p\in\mathcal{H}$.
\end{thm}

The existence is a direct application of the first named author's
result~\cite{ChenC11}. However, we have to assume that the total space of the test configuration is smooth.  We hope the result can be extended to singular test configurations accordingly.

In\cite{Chenbounded}, the first author took another approach to construct the geodesic ray. Using techniques in\cite{Chenbounded},  the smooth condition here can be reduced to the lower bound of the Riemann curvature of the total space. 

\proof Consider a smooth test configuration over a disc:
$(\mathcal{L}\rightarrow\mathcal{M}\rightarrow
D)\hookrightarrow(\mathcal{O}(1)\rightarrow P^N\times D\rightarrow
D)$. Assume the total space is smooth. i.e, $\mathcal{M}\subset
P^N\times D$ is smooth. Let $\Omega$ be the Fubini-study metric on
$P^N\times D$. Actually, it means the pull back of Fubini-study
metric on $P^N$ by projection: $P^N\times D\rightarrow P^N$.

Now solve the equation
\begin{eqnarray}\label{Basic}
(\Omega+\sqrt{-1}\partial\bar\partial\psi)^{n+1} &=& 0 \texttt{
     on    }\mathcal{M}
\\ \psi &=&0\texttt{    on    }\partial \mathcal{M}
\end{eqnarray}
According to~\cite{ChenC11}, this equation has a $C^{1,1}$ solution(
It is not exactly the same situation as in~\cite{ChenC11}, but the techniques are the same). The following
shows that: This solution corresponds to a geodesic ray in the K\"ahler
class $c_1(L)$.

The $C^*$ action on $\mathcal{M}$ induces a biholomorphic map
$i:(L_1,M_1)\times (D-0)\rightarrow (\mathcal{L},\mathcal{M})-M_0$. Now
$i$ maps $(e,x,z)\in (L_1,M_1)\times (D-0)$ to
$z\circ(e,x,1)\subset (\mathcal{L},\mathcal{M})$. $z\circ$ is the $C^*$ action of test configuration, and $(e,x,1)\in(L_1,M_1)$. The map $i$ pulls the equation to
\begin{equation}
(i^*\Omega+\sqrt{-1}\partial\bar\partial i^*\psi)^{n+1}=0
\end{equation}
on $M_1\times (D-0)$, with boundary condition $i^*\psi=0$ on
$M_1\times S^1$.

Let $\omega=\Omega|_{M_1}$, and $\pi:M_1\times (D-0)\rightarrow M_1$ be
the projection, then
\begin{prop} $i^*\Omega=\pi^*\omega+\sqrt{-1}\partial\bar\partial\eta$ for some smooth function $\eta$.
\end{prop}
\proof Let $h$ be the Fubini-Study hermitian metric on
$\mathcal{O}(1)\rightarrow P^N$. So
$\Omega=-\sqrt{-1}\partial\bar\partial\log h$ and
$i^*\Omega=-\sqrt{-1}\partial\bar\partial\log i^*h$. Note
$\pi^*\omega=-\sqrt{-1}\partial\bar\partial\log h_1$. $h_1$
is the pull back of the hermitian metric on line bundle
$L_1\rightarrow M_1$ by trivial projection $\pi:(L_1,M_1)\times
(D-0)\rightarrow (L_1,M_1)$. So
$i^*\Omega=\pi^*\omega+\sqrt{-1}\partial\bar\partial\log\frac{h_1}{i^*h}$
and $\eta=\log\frac{h_1}{i^*h}$. $\Box$

\begin{prop} $\varphi=\eta+i^*\psi$ is a geodesic ray.
\end{prop}
\proof We have showed
$(\pi^*\omega+\sqrt{-1}\partial\bar\partial\varphi)^{n+1}=0$ on
$M\times (D-0)$. It remains to show the $S^1$
invariance of $\varphi$. First, we check the $S^1$ invariance of $\eta$. By assumption, $S^1$ action on $\mathcal{O}(1)\rightarrow P^N\times C$ is unitary. So the $h$ is preserved by $S^1$ action. This immediately implies $\eta=\log\frac{h_1}{i^*h}$ is $S^1$ invariant. Now we check $\psi$. $\psi$ is
$S^1$ invariant because the boundary condition $\psi=0$ is $S^1$ invariant, and
the uniqueness of Monge Ampere solution. In another word, for the unique solution, the $S^1$
symmetric on the boundary will force the $S^1$ symmetry in the
interior. Now both $\eta$ and $\psi$ are $S^1$ invariant, so is $\varphi$. $\Box$

So far, we have associated a relative $C^{1,1}$ geodesic ray to the test configuration. 
The ray starts from a fixed point $p$, because we solved the equation with boundary condition $\psi=0$. However, for another arbitrary point $q$, one can go back to the equation~\ref{Basic}, solve $\psi=\psi_0$ on $\partial\mathcal{M}$ and obtain the relative $C^{1,1}$ ray from $q$. $\psi_0$ is the $S^1$ extension of the potential difference between $q$ and $p$. 
$\Box$

In~\cite{Tianarezzo}, Arezzo and Tian
constructed an analytic geodesic ray from a
test configuration when the central fiber is analytic. Such test configurations in~\cite{Tianarezzo} are simple test configurations~\ref{SimpleTest}. Using the openness theorem~\ref{Open}, we know that there are smooth geodesic rays near the ray they constructed.\\

Back to the question: given a geodesic ray, how to construct a test configuration which represents the same degeneration? Donaldson's construction of toric degenerations~\cite{DonaldsonToric} is very inspiring: He chose piece wise linear functions to approximate an arbitrary direction, and the piece wise linear function leads to a well defined test configuration. In principle, one can think the degenerations represented by test configuration are dense in all possible geometrical degenerations. Donaldson's construction is a method to choose good approximation, which reflects the same character of degeneration.   

\subsection{Special cases: geodesic line and Toric variety}
One example of geodesic ray is the geodesic line generated by a holomorphic vector field.  Let $M$ be K\"ahler manifold with K\"ahler form $\omega$. Let $X$ be a holomorphic vector field such that: $X=f^{,\alpha}\frac{\partial}{\partial w^\alpha}$ for some real potential $f$ and $Im(X)$ is killing vector field. Let $\sigma(t)$ be the flow generated by $Re(X)=\nabla_\omega f$. Then the 1-parameter family $\omega_{\rho(t)}=\sigma(t)^*\omega$ is a geodesic line, $t\in(-\infty,\infty)$. 

Another special case is when the manifold is a toric variety. For a toric variety, there is an associated polytope. In detail, there is biholomorphic map $f:M^\circ=C^n/2\pi i Z^n\rightarrow P^\circ\times T^n$. $M^\circ$ is an open dense subset of $\mathcal{M}$ where the toric action is free. $P$ is a polytope in $R^n$ satisfying Delzant condition. If we write a Toric-invariant K\"ahler metric $\omega|_{M^\circ}=i\partial\bar\partial f$, then there is a map from $C^n/2\pi i Z^n$ to $P^\circ\times T^n$: $(u,v)\rightarrow (x=\frac{\partial f}{\partial u},y=v)$. Under this map , the K\"ahler form $\omega$ is translated into $dx\wedge dy$. The complex structure is translated into

\begin{equation}
J=\left(
\begin{array}{ccc}
  0   &G   \\
  G^{-1}   &0 
\end{array}
\right)
\end{equation}
$G_{ij}=\frac{\partial^2 g}{\partial x_i\partial x_j}$,  $g(x)+f(u)=\sum x_i u_i$ at $x=\frac{\partial f}{\partial u}$. In another word, in the symplectic chart, the complex structure has a potential $g$. 

This transformation is really helpful for the geodesic equation. The geodesic equation, in the polytope representation, is linear for Complex structure potential $g(t)$. i.e,
\begin{equation}
\ddot g(t)=0
\end{equation}
This immediately implies the existence of smooth geodesics connecting any two toric metrics. It is just the linear interpolation of the two end potentials.

\section{Connection between algebraic notions and geometric notions}
\subsection{Algebraic ray and geodesic ray}
Test configurations can be viewed as algebraic rays. The induced geodesic rays are parallel to the algebraic ray. 

\begin{defn}Two rays $\rho_1(t)$ and $\rho_2(t)$ in the space of K\"ahler metrics are called parallel if $\rho_1(t)-\rho_2(t)$ is uniformly bounded.
\end{defn}
The equality $\varphi=\eta+i^*\psi$ can be interpreted
geometrically. $\eta$ represents the degeneration of the metric
from the algebraic $C^*$ action. $\psi$ is the difference between
the algebraic ray and the differential geometric ray. Notice that $\psi$
is $C^{1,1}$ bounded. We will elaborate above statement in the following:

Recall that $(L,M)=(L_1,M_1)\hookrightarrow (\mathcal{O}(1),P^N)$ is
embedding. The group $GL(N+1,C)$ acts on $(\mathcal{O}(1),P^N)$. If
one looks at the dual bundle of $\mathcal{O}(1)$ (i.e. the
universal bundle $\{(e,x)\in C^{N+1}\times P^N: e=\lambda x\}$), the
action is simply $A(e,x)=(Ae,Ax), A\in GL(N+1,C)$. The natural dual
map between $\mathcal{O}(1)$ and universal bundle passes the action
from one to the other.

Consequently, the action acts on the hermitian metric of
$\mathcal{O}(1)$, thus on its curvature. The following lemma shows
it preserves the positivity of the hermitian curvature.

\begin{lem}Let $A\in GL(N+1,C)$, $h$ is the Fubini-Study hermitian metric on $\mathcal{O}(1)$, then
$-i\partial\bar\partial\log A^*h>0$
\end{lem}
\proof It suffices to prove that the action preserves the
negativity of curvature on the universal bundle. Under the action
$A$, the metric of $e=(X_0,X_1,...,X_N)\in \mathcal{O}(-1)$
changes into $||Ae||^2$ from standard Fubini-Study metric
$||e||^2$. Notice that the action $A^{-1}UA$ for $U\in U(N+1)$ is
transitive on $P^N$ and preserves the $A^*h$. So one just needs to
show negativity at one point. Lets consider the point
$p=A^{-1}(1,0,...,0)^t$, and
$e=(X_0,...,X_{i-1},1,X_{i+1}...,X_N)$. At the point $p$,
\begin{equation}
-\sqrt{-1}\partial\bar\partial\log
||Ae||^2=-\sqrt{-1}\sum_{j=1}^n\sum_{k,l\neq i} A_{jk}\bar
A_{jl}dX_k\wedge d\bar X_l
\end{equation}
To show the positivity, it suffices to show that the null space of the
matrix $A_{jk}, j\neq 1,k\neq i$ must be empty. If 
$v=(\alpha_0,...,\alpha_{i-1},\alpha_{i+1},...\alpha_N)$ is a
null vector, then the
vector
$Av^t$ must
be of form $(c\neq 0,0,0,...,0)$, because of
non-singularity of $A$. By scaling $c=1$, $A$ will map two vectors
to $(1,0,...,0)$, contradiction. $\Box$

The consequence is: $GL(N+1,C)$ action induces a finite
dimensional subspace $\mathcal{H}_N\subset\mathcal{H}$.
$\mathcal{H}_N$ consists of those metrics obtained by the
$GL(N+1,C)$ action.

The space $\mathcal{H}_N$ is a symmetric space. Its dual is the
unitary group $U(N+1)$. Under the natural metric of symmetric spaces, the $C^*$ action (as a 1-parameter
family of metrics) is a geodesic ray in $\mathcal{H}_N$. It is
interesting to consider the limit of these algebraic rays when
one raise the dimension of ambient space $P^N$(we can raise the power $k$ of $\mathcal{L}^k$ and do Kodaira embedding, then pull the ray back to the class $c_1(L)$ by dividing out the scalar $k$). First, it is easy to derive that all the embedding induce the same geometric geodesic ray.
\begin{lem} Different embedding of a test configuration into projective spaces induce the same geodesic ray provided the rays start at the same point.
\end{lem}
\proof By different embedding, one essentially raise the power $k$ of $\mathcal{L}^k\rightarrow\mathcal{M}\rightarrow D$ and then use sections of $H^0(\mathcal{M},\mathcal{L}^k)$ to embed $\mathcal{L}^k\rightarrow\mathcal{M}$ into $\mathcal{O}(1)\rightarrow P^N\times D$. The Fubini-Study metric naturally induces a metric on $\mathcal{L}^k$, which has curvature in class $kc_1(\mathcal{L})$. To get a geodesic ray in the K\"ahler class $c_1(L)$, one takes the $k$-th root of the Fubini metric on $\mathcal{L}^k$ to get a hermitian metric $h_k$ on $\mathcal{L}$. Notice that $\log\frac{h_k}{h_n}$ is the potential difference of the background metric $\Omega_k$ and $\Omega_n$. When we solve the Monge Ampere equation, this finite difference goes into the $C^{1,1}$ solutions $\phi_k$ and $\phi_n$. Thus the ray potential $\eta_k+i^*\phi_k=\eta_n+i^*\phi_n$. $\Box$

As $k\rightarrow\infty$, it is expected these algebraic rays should converge to the geometric geodesic ray. This is a natural extension of the classical problem: Use Bergman metrics to approximate a given K\"ahler metric. There is extensive literature on this topic, c.f. Tian\cite{BergmanTian}, Zelditch\cite{Zelditch}, Lu\cite{Lu}, Song\cite{SongZelditch}.

\subsection{Bounded ambient geometry and test configuration}
In~\cite{Chenbounded}, the first named author defined bounded ambient geometry to study geodesic rays. Briefly speaking, a geodesic ray is called to have bounded ambient geometry if the following holds: There exists a metric $\tilde g$ on $M\times S^1\times[0,\infty)$ such that the ray has a $C^{1,1}$ relative potential under $\tilde g$, and $\tilde g$ has uniformly bounded curvature.

The geodesic ray induced by a smooth test configuration always has bounded ambient geometry. To see this, one restrict the metric $\Omega+idz\wedge d\bar z$ to the punched part $\mathcal{M}-M_0$. Since $\Omega+idz\wedge d\bar z$ has bounded geometry on $\mathcal{M}$, the restriction clearly has bounded geometry. The punched part is holomorphic identified with $M\times S^1\times[0,\infty)$, thus the ray has bounded ambient geometry.
Actually, it is a stronger version of bounded ambient geometry since the metric $\tilde g$ on $M\times S^1\times [0,\infty)$ can be compactified into disc fiberation. In general cases of bounded ambient geometry, this is not necessarily true.

In~\cite{Chenbounded}, it is proved that: Let $\rho(t)$ be a geodesic ray with bounded ambient geometry, then for any other potential $\phi_0$, there is a unique relative $C^{1,1}$ geodesic ray starting from $\phi_0$ and parallel to $\rho(t)$. Alternatively, we can use this to derive the existence of geodesic rays, based on the algebraic ray.

\subsection{Futaki invariant, $\yen$ invariant and stability}
Futaki invariant is introduced by Futaki, on the manifold with positive chern class $c_1>0$. Later, Calabi extended the definition to general K\"ahler manifold. Ding and Tian generalized Futaki invariant for a class of singular varieties~\cite{DTFutaki} and Donaldson defined Futaki invariant for test configurations. 

The classical definition of Futaki invariant is the following:
Let $M$ be a K\"ahler manifold with K\"ahler metric $\omega$. Let $X$ be a holomorphic vector field on $M$. Let $h$ be the solution of $\Delta h=R-\underline{R}$. Futaki invariant is a linear functional: $\mathcal{F}(X)=\int_M X(h)\omega^n$.  The definition is independent with the metric $\omega$ chosen in a fixed class. In particular, when $X=f^{,\alpha}\frac{\partial}{\partial w^\alpha}$, $\mathcal{F}(X)=\int_M f^{,\alpha}h_{,\alpha}\omega^n=\int_M f(\underline{R}-R)\omega^n$.

Ding and Tian~\cite{DTFutaki} generalized the Futaki invariant to a class of singular varieties. Briefly speaking, they embed the variety into a projective space $P^N$, and consider the restriction of ambient holomorphic vector fields tangent to the variety on regular points. Also they consider the restriction of ambient Fubini-study metric $\omega$ and define Futaki invariant in similar fashion. 

In test configuration, Donaldson's algebraic definition of Futaki invariant is:
Let $\mathcal{L}\rightarrow\mathcal{M}\rightarrow D$ be a test configuration. Consider the $C^*$ action on the central fiber $L_0\rightarrow M_0$, and its powers $L_0^k\rightarrow M_0$. Let $d_k=\dim H_k=\dim H^0(M_0;L_0^k)$ and $w_k$ be the weight of the $C^*$ action on highest exterior power of $H_k$. Then $F(k)=w_k/kd_k$ has an expansion
\begin{equation}
F(k)=F_0+F_1k^{-1}+F_2k^{-2}+......
\end{equation}
The coefficient $F_1$ is called the Futaki invariant of the $C^*$ action on $(L_0,M_0)$.
He proved that if the central fiber is smooth, then the algebraic Futaki invariant agrees with the classical Futaki invariant. 

Using Futaki invariant, Donaldson defined stability. A pair $(L,M)$ is K-stable if: For each test configuration for $(L,M)$ (i.e, $(L_1,M_1)=(L,M)$), the Futaki invariant of the $C^*$ action on $(L_0,M_0)$ is less than or equal to zero, and the equality only occurs when the configuration is a product configuration.

This algebraic definition agrees with an early geometric definition of K-stability by Ding and Tian. In~\cite{DTFutaki}, they used a $C^*$ action of $P^N$ to obtain the limit of the varieties $M_t$, then studied the Futaki invariant of the limiting variety $M_0$. The spirit is similar to Donaldson's setup of test configuration.

Notice that in test configuration, the stability is to check the Futaki invariant of the central fiber. However, one would like to have some criterion that doesn't need a specific central fiber. Just as the bounded ambient geometry only concerns behavior before reaching the limit, the $\yen$ invariant is a nice notion parallel to Futaki invariant and doesn't need a specific central fiber. Following~\cite{Chenbounded}.
\begin{defn} For a smooth geodesic ray $\rho(t)$, $\yen$ invariant is defined to be
\begin{equation}
\yen=\lim_{t\rightarrow\infty}\frac{dE}{dt}=\lim_{t\rightarrow\infty}\int\frac{\partial\rho}{\partial t}(\underline{R}-R)\omega_\rho^n
\end{equation}
\end{defn}

The K-engery is convex along geodesics. So $\frac{dE}{dt}$ is monotone and the limit exists.

The first named author defined the notion of geodesic stability by $\yen$ invariant:  $M$ is weakly geodesically stable if every geodesic ray has nonnegative $\yen$ invariant. This is parallel to K-stability in test configurations. However, geodesic rays represent all possible geometrical degenerations. So it is possible that the geodesic ray might detect some instabilities which test configuration cant detect.

To clarify the analogy, we prove that: In the case of simple test configurations, the $\yen$ invariant agrees with the Futaki invariant. 
\begin{thm}For simple test configuration, if the geodesic ray is smooth regular, then $\yen$ invariant agrees with Futaki invariant\footnote{It is the same up to a sign}.
\end{thm}
\proof By definition of simple test configuration, the central fiber is smooth. Following~\cite{DonaldsonToric}, the algebraic Futaki invariant is exactly the classical Futaki-invariant applying to the $C^*$ action holomorphic vector field in the central fiber.

Let $\omega_c$ be the restriction of $\Omega+i\partial\bar\partial\phi$ on $M_0$. The $S^1$ action of the $C^*$ action is a hamiltonian action on $M_0$. Let $f$ be the hamiltonian. In another word, $df=i_v\omega_c$, where $v$ is the $S^1$ action vector field.  The Futaki-invariant of the $C^*$ action is 

\begin{equation}
\nu=\int f(\underline{R}-R)\omega_c^n
\end{equation}

Now we look at $\yen=\lim_{t\rightarrow\infty}\int\frac{\partial\rho}{\partial t}(\underline{R}-R)\omega_\rho^n$. If we apply a diffeomorphism to each $M$ in $M\times[0,\infty)\times S^1$, i.e, identify the $M\times[0,\infty)\times S^1$ with $\mathcal{M}-M_0$, then
\begin{equation}
\lim_{t\rightarrow\infty}\omega_\rho=\omega_c, \lim_{t\rightarrow\infty}R_\rho=R_{\omega_c}
\end{equation}

So it suffices to show $\lim_{t\rightarrow\infty}\frac{\partial\rho}{\partial t}=-f$. 

Notice the following fact: In $M\times[0,\infty)\times S^1$, the solution foliation induces an $S^1$ action, which is moving along the leaf in $S^1$ direction. By identifying the fiber $M_t$ with $M_{t\theta}$ where $|\theta|=1$, the $S^1$ action is hamiltonian action with hamiltonian $\frac{\partial\rho}{\partial t}$, under the symplectic form $\omega_\rho$. By translating this into the context of $\mathcal{M}$, we have: If we identify the fiber $M_t$ with $M_{t\theta}$ via the $S^1$ action of the $C^*$ action, then the $S^1$ action induced by foliation is hamiltonian action with hamiltonian $\frac{\partial\rho}{\partial t}$, under symplectic form $\omega_\rho$.  Now we take limit of the identification towards the central fiber, the $S^1$ action induced by foliation converges to the $S^1$ action of the $C^*$ action on the central fiber. In the picture of test configuration, the limit of the $S^1$ action induced by foliation is trivial in the central fiber. However, because of the distortion created by the identification, the limit under this identification is  the reverse of $S^1$ action of the $C^*$ action on central fiber. 

Therefore, the limit of the hamiltonian $\frac{\partial\rho}{\partial t}$ is the hamiltonian of the limiting action. So $\lim_{t\rightarrow\infty}\frac{\partial\rho}{\partial t}=-f$ and the theorem is proved. $\Box$

It is a well known conjecture that existence of constant scalar curvature metric or extremal metric is related to the stability of K\"ahler manifolds. Yau pointed out this in 1980s. From then on, there has been much progress in this topic. Interested readers may consult the rich literature in this area, Tian\cite{inventionpaper97},  Donaldson~\cite{Donaldsonlowerbound} , Mabuchi\cite{Stability}, Paul-Tian\cite{PaulTian},  chen-tian\cite{ChenTian}, Zhou-Zhu\cite{Zhu}...

\section{Monge Ampere equation on Simple test configurations}
Following Donaldson's idea~\cite{DonaldsonSetup}, this section
extends the correspondence in~\cite{DonaldsonSetup} to the case of
Monge Ampere equation on simple test configurations.

But to explain the background and the motive, we start with a review
on Donaldson's result. $M$ is a K\"ahler manifold with a given K\"ahler
form $\omega$. We solve the equation
$(\pi^*\omega+i\partial\bar\partial\phi)^{n+1}=0$ on $M\times D$
with boundary condition $\phi=\phi_0$ on $M\times\partial D$. $\pi$
is the natural projection to $M$.

Donaldson and Semmes independently constructed the following
manifold $W\rightarrow M$. $W$ is glued by local holomorphic
cotangent bundle over $M$. There exists a lifting of $M$ into $W$
for every K\"ahler metric $\omega+i\partial\bar\partial\phi$. If one
take the lifting of $M\times D$ into $W\times D$ by the solution
$\omega+i\partial\bar\partial\phi$, then one will obtain a family of
holomorphic discs. These discs are the lifting of the foliation induced by
the degenerated form $\pi^*\omega+i\partial\bar\partial\phi$. Vice
versa, if one has the family, then it can induce a solution to Monge Ampere equation. This correspondence is useful. It relates the PDE regularity to the moduli
space of holomorphic discs.

However, Donaldson's construction only works for a product like $M\times
D$. But a test configuration of real interest is not a product space. So the previous construction won't work directly. We solve this problem by taking a new point of view on the old construction: View $W\times D$ as a global construction over $M\times D$. Then we can derive an analogy in non-product case.

\subsection{Construction of $\mathcal{W}\rightarrow\mathcal{M}$}
Recall a test configuration is simple (definition~\ref{SimpleTest}) if: The total space $\mathcal{M}$ is smooth ($\mathcal{M}$ is a smooth sub-manifold of $P^N\times D$) and the
projection $\pi:\mathcal{M}\rightarrow D$ is submersion everywhere.

From the definition, simple test configuration is fiberation over
disc. Each fiber is smooth because $\pi:\mathcal{M}\rightarrow D$ is
submersion everywhere.

Let $\mathcal{M}$ be a simple test configuration. 
We solve $(\Omega+i\partial\bar\partial\phi)^{n+1}=0$ on $\mathcal{M}$. Since $\pi:\mathcal{M}\rightarrow D$ is submersion everywhere, so
$\mathcal{M}$ is locally product space. To see this explicitly in the complex
coordinates: First, choose a complex
coordinate $\{x_0,...,x_n\}$ for $U\subset\mathcal{M}$. The
projection $z=z(x_0,...,x_n)$ is holomorphic and
$\frac{\partial z}{\partial x_i}\neq 0$ by assumption of submersion. Now one
can easily cook up a tuple $\{z,x_{i_1},...,x_{i_n}\}$ such that the
transition between $\{z,x_{i_1},...,x_{i_n}\}$ and $\{x_0,...,x_n\}$ is non-degenerate.
$\{z,x_{i_1},...,x_{i_n}\}$ is the product holomorphic coordinate we are looking for.

In the future, such product coordinate is denoted by $(z,w)$ with
$z\in D$ and $w\in M_z$. Cover $\mathcal{M}$
with local product charts $U_i$. On $U_i$, suppose the
$\Omega=i\partial\bar\partial\rho_i$. Write $T^*\mathcal{M}/T^*C$
over $U_\alpha$ by local coordinates $(z,w,q)$. We glue these charts together, and define the transition
between $(z,w,q)$ over $U_\alpha$ and $(v,x,p)$ over $U_\beta$:
\begin{eqnarray}
\nonumber z&=&v\\
\nonumber w&=&w(v,x)\texttt{  as defined in   }\mathcal{M}\\
q_j&=&p_i\frac{\partial x_i}{\partial
w_j}+\frac{\partial(\rho_\beta-\rho_\alpha)}{\partial w_j}
\end{eqnarray}
One can verify these local charts
$(z,w,q)$ glue up to a complex manifold
$\mathcal{W}\rightarrow\mathcal{M}$. Define a form
$\Theta$ on each fiber of $\mathcal{W}\rightarrow D$,
\begin{equation}
\Theta|_{W_t}=dq_i\wedge dw_i
\end{equation}
$\Theta$ is well defined only on the fiber, so $\Theta|_{W_t}$ is a family of forms.

The real part of $\Theta$ is a symplectic form on $W_t$. So $W_t$ is a symplectic manifold and we can talk about Lagrangian sub-manifolds of $W_t$.

\begin{defn}For a
Lagrangian sub-manifold $L_t$, $L_t$ is called LS-submanifold if $\Theta|_{L_t}$ is
non-degenerate. $L_t$ is called LS-graph if it is LS-submanifold and also be a graph over $M_t$.\end{defn}
By straightforward calculation, one can see LS-graphs are of forms
$\partial\phi$ for some real potential $\phi$ on $M_t$, and
$\Theta|_{L_t}=\partial\bar\partial\phi$.

Our main result in section 4 is:
\begin{thm}Let $\mathcal{M}$ be a simple test configuration. There
is an associated manifold $\mathcal{W}\rightarrow\mathcal{M}$ such
that:
\begin{enumerate}
  \item A smooth solution $\phi$ of
  $(\Omega+i\partial\bar\partial\phi)^{n+1}=0, \phi=\phi_0$ on
  $\partial\mathcal{M}$ induces a family of holomorphic discs
  $G:M\times D\rightarrow\mathcal{M}\rightarrow\mathcal{W}$ factoring through the foliation on $\mathcal{M}$,
  such that the image of $G(\cdot,z)$ is a LS-graph in $W_z\rightarrow M_z$ for all $z$ and
  $\bigcup_{z\in\partial D} G(\cdot,z)$ is a totally real
  sub-manifold of $\mathcal{W}$.
  \item If a family of holomorphic discs $G:M\times
  D\rightarrow\mathcal{W}$ respects the projection $\mathcal{W}\rightarrow D$, i.e, $\pi\circ G:M\times D\rightarrow D$ is a projection to
  $D$. Also assume it satisfies the boundary condition $
  G(\cdot,z)=\Lambda_{z,\phi_0}$ for $z\in\partial D$, where
  $\Lambda_{z,\phi_0}$ is the lifting of $M_z$ by metric
  $\Omega+i\partial\bar\partial\phi_0$, then the image of
  $G(\cdot,z)$ is a LS-submanifold in $W_z$ for all $z$. Moreover, if assuming
  these images are LS-graphs, then the family projects to a foliation of $\mathcal{M}$ and
  induces a smooth solution $\phi$ to $(\Omega+i\partial\bar\partial\phi)^{n+1}=0$
  with $\phi=\phi_0$ on $\partial\mathcal{M}$.

\end{enumerate}
\end{thm}

\subsection{One side of the Correspondence}
Now suppose there is a smooth solution $\phi$ for
$(\Omega+i\partial\bar\partial\phi)^{n+1}=0$ on $\mathcal{M},
\phi=\phi_0$ on $\partial\mathcal{M}$, with
$\Omega+i\partial\bar\partial\phi$ positive on $M_t$.

In local product coordinates $(z,w)$ of $\mathcal{M}$, write
$\Omega+i\partial\bar\partial\phi=i\partial\bar\partial f$. Since
$\Omega+i\partial\bar\partial\phi$ has rank $n$, it has a 1-complex
dimension kernel. Let $X=\frac{\partial}{\partial
z}+\eta^\alpha\frac{\partial}{\partial w^\alpha}$ be in kernel of
$i\partial\bar\partial f$, then
\begin{equation}
0=\partial\bar\partial f(\frac{\partial}{\partial
z}+\eta^\alpha\frac{\partial}{\partial w^\alpha})=(\eta^\alpha
f_{\alpha\bar\beta}+f_{z\bar\beta})dw^{\bar\beta}+(\eta^\alpha
f_{\alpha\bar z}+f_{z\bar z})d\bar z
\end{equation}
so
\begin{eqnarray}
\eta^\alpha&=&-f_{z\bar\beta}f^{\alpha\bar\beta}\\
 f_{z\bar
z}&=&-\eta^\alpha f_{\alpha\bar z}
\end{eqnarray}
Now, direct calculation shows
\begin{equation}
[X,\bar X]=(\frac{\partial\eta^{\bar\beta}}{\partial
z}+\eta^\alpha\frac{\partial\eta^{\bar\beta}}{\partial
w^\alpha})\frac{\partial}{\partial
w^{\bar\beta}}-(\frac{\partial\eta^{\alpha}}{\partial\bar
z}+\eta^{\bar\beta}\frac{\partial\eta^\alpha}{\partial
w^{\bar\beta}})\frac{\partial}{\partial w^\alpha}=0
\end{equation}
this means the kernel distribution is holomorphicly parametrized by
$z\in D$. So a smooth solution implies a foliation of $\mathcal{M}$
by holomorphic discs.

The $\mathcal{M}$ can be lifted to a graph in $\mathcal{W}$, using
the form $\Omega+i\partial\bar\partial\phi$. In detail, on local
product charts $U_i$, $\Omega=i\partial\bar\partial\rho_i$, we can lift
$\mathcal{M}$ to graph $\partial(\rho_i+\phi)$ in each fiber. The
lift is well defined globally due to the way we glue
$\mathcal{W}$.

In~\cite{DonaldsonSetup}, Donaldson showed in the lifting of
$\mathcal{M}$, the foliation is lifted up to a family of holomorphic
discs in $\mathcal{W}$, and these holomorphic discs take boundary value in a totally real
sub-manifold $\Lambda_{\phi_0}$. The same technique works here.

In summary,
\begin{thm}For a simple test configuration, the smooth solution of
the Monge ampere equation induces a foliation of holomorphic discs
on $\mathcal{M}$ which lift up to a family of holomorphic discs with
in $\mathcal{W}$. These discs have boundary in a totally real
sub-manifold.
\end{thm}
\proof As above. $\Box$
\subsection{The other side of the correspondence}
It is reasonable to consider the reverse correspondence locally.
We
have the following theorem:

\begin{thm}Suppose $G:D\times U\rightarrow\mathcal{W}$ is a smooth map which respects
the projection and holomorphic in $D$. Assume for all
$\tau\in\partial D$, $U$ is mapped to be LS-graph and this LS-graph family has a global
potential $\phi_0$. Then for each
$\tau\in D$, $G$ maps $U$ to an immersed LS-submanifold in
$\mathcal{W}$. Moreover, if assuming these LS-submanifolds are LS-graphs, then this
family induces a smooth solution to the Monge Ampere equation with boundary condition $\phi=\phi_0$.
\end{thm}
In above theorem, $U$ is an open set of real dimension $2n$.  $G:D\times
U\rightarrow \mathcal{W}$ is smooth and respects the
projection. In another word, for $\pi:\mathcal{W}\rightarrow D$, $\pi\circ G$
is identity on $D$. $G$ is holomorphic in $D$
variable. For each $\tau\in\partial D$, $U$ is mapped to be a
LS graph over $M_\tau$ and this LS-graph family have a global
potential $\phi_0$. This just means these LS-graphs are lifting of $\mathcal{M}$
using $\Omega+i\partial\bar\partial\phi_0$ on the boundary. 

\proof Consider $G^*\Theta$ on $D\times U$. $\Theta$ is well defined
on fibers $W_t$, so $G^*\Theta$ is well defined on fibers $U_t$
in $D\times U$. We should view $G^*\Theta$ as a family of forms on $U_t$.
Denote real coordinates on $U$ by $q_j$, write
$G^*\Theta=(r_{jk}+is_{jk})dq_j\wedge dq_k$. It is straightforward
to show $r_{jk}+is_{jk}$ is holomorphic function over $D$: Let
$(z,q)$ be coordinates on $D\times U$. Let $(v,x,p)$ be a local
coordinates in $\mathcal{W}$. The map $G$ is $v=z, x=x(z,q),
p=p(z,q)$. $G$ is holomorphic, so $\frac{\partial x}{\partial \bar
z}=0,\frac{\partial p}{\partial\bar z}=0$. Now
$\Theta|_{W_t}=dp_i\wedge dx_i, G^*\Theta|_{U_t}=\frac{\partial
p_i}{\partial q_j}\frac{\partial x_i}{\partial q_k}dq_j\wedge dq_k$,
therefore $\frac{\partial}{\partial \bar
z}(r_{jk}+is_{jk})=\frac{\partial}{\partial \bar z}\frac{\partial
p_i}{\partial q_j}\frac{\partial x_i}{\partial q_k}=0$.

On the boundary $\tau\in\partial D$, $G$ maps $U$ to LS-graphs. But
$\Theta$ is purely imaginary on LS-graphs, so $G^*\Theta$ is also
purely imaginary. A holomorphic function on the disc with pure
imaginary value on $\partial D$ must be constant, so
$r_{ij}+is_{ij}$ must be constant on every disc in $D\times U$. This
also implies the Jacobi of the map $G(\tau,\cdot):U\rightarrow
W_\tau$ is non-degenerate, since the pull back image $G^*\Theta$ is
non-degenerate. It follows that the image $G(\tau,U)$ is an immersed
LS-submanifold.

Now assume $G(\tau,U)$ is actually a LS-graph, i.e, the projection
$\pi\circ G(\tau,\cdot)$ is diffeomorphism. Following~\cite{ChenTian}, we find a global
potential for this family of LS-graphs (modulo the local potential
of the background metric).

First, consider the case when $U$ is a very small open ball. Let
$D_\alpha$ be a small open set in $D$. Without loss of generality,
$G$ maps $D_\alpha\times U$ into a single chart in $\mathcal{W}$.
Since they are LS-graphs, one can solve a real potential
$\varphi_\alpha$ for this family in the local product chart by
$\frac{\partial\varphi_\alpha}{\partial x_i}=p_i$.
$\varphi_\alpha$ is unique up to a smooth function in $z\in D$.

Choose a finite covering $D_\alpha\subset D$, and make $U$ so
small such that $D_\alpha\times U$ all fit in single charts in
$\mathcal{W}$. This can be done if one fixes a finite chart
covering of $\mathcal{W}\rightarrow D$ in first place and then
replace $U$ by small subset if necessary. Solve the potential
$\varphi_\alpha$ respectively in each $D_\alpha\times U$, and the
geometry of $\mathcal{W}$ implies
$\partial(\varphi_\alpha-\rho_\alpha)=\partial(\varphi_\beta-\rho_\beta)$
on every fiber $M_t$ of $\mathcal{M}$. So on each fiber, the
difference
$(\varphi_\alpha-\rho_\alpha)-(\varphi_\beta-\rho_\beta)$ must be
constant. It follows that $\varphi_\alpha-\rho_\alpha$ differ with
$\varphi_\beta-\rho_\beta$ by a smooth real function of $z$ on
intersection. The fact $H^1(D,\mathcal{S})=0$, ($\mathcal{S}$ is
the sheaf of $C^\infty$ functions) implies one can adjust
$\varphi_\alpha$ by function of $z$ such that
$\varphi_\alpha-\rho_\alpha=\varphi_\beta-\rho_\beta$. Therefore
they give the global potential $\phi=\varphi_\alpha-\rho_\alpha$. $\phi$ is unique up to a function of $z$ on $D$.

The next step is to make $\phi$ satisfy the boundary condition
$\phi=\phi_0$. Let $X=\frac{\partial}{\partial
z}+\eta^\alpha\frac{\partial}{\partial w^\alpha}$ be the tangential
vector of the foliation $\pi\circ G:D\times
U\rightarrow\mathcal{M}$. There exists a 1-1 form $\Omega'$ on
$\mathcal{M}$ such that $i_X\Omega'=0$ and its restriction to $M_t$
is $i\partial\bar\partial\varphi_\alpha|_{M_t}=\Theta|_{L_t}$.
Locally, $\Omega'=i(\frac{\partial^2\varphi}{\partial w_\alpha
w_{\bar\beta}}dw_\alpha dw_{\bar\beta}+\zeta^\alpha dw_\alpha d\bar
z+\zeta^{\bar\beta}dw_{\bar\beta}dz+hdzd\bar z)$, where
$\zeta^\alpha=-\eta^{\bar\beta}\varphi_{\alpha\bar\beta}$ and
$h=\eta^\alpha\eta^{\bar\beta}\varphi_{\alpha\bar\beta}$.

Let $(v,q)$ be coordinates on $D\times U$, $q$ as real coordinates.
$(z,w)$ are local coordinates on $\mathcal{M}$. We have 
$\eta^{\bar\beta}=\frac{\partial w^{\bar\beta}}{\partial\bar v}$.
Let $\rho$ be local potential for background metric $\Omega$, and
$\varphi=\rho+\phi$. The disc family in $\mathcal{W}$ is holomorphic
implies $\frac{\partial}{\partial \bar
v}\frac{\partial\varphi}{\partial w_\alpha}=0$, therefore
\begin{equation}
0=\frac{\partial}{\partial \bar v}\frac{\partial\varphi}{\partial
w_\alpha}=\frac{\partial^2\varphi}{\partial w_\alpha \partial\bar
z}+\frac{\partial^2\varphi}{\partial w_\alpha\partial
w_{\bar\beta}}\eta^{\bar\beta}
\end{equation}
So $\zeta^\alpha=\frac{\partial^2\varphi}{\partial w_\alpha
\partial\bar z}$,
$\Omega'=i(\partial\bar\partial\varphi+(h-\varphi_{z\bar z})dzd\bar
z)=i(\partial\bar\partial(\rho+\phi)+(h-\rho_{z\bar z}-\phi_{z\bar
z})dzd\bar z)=\Omega+i\partial\bar\partial\phi+i(h-\rho_{z\bar
z}-\phi_{z\bar z})dzd\bar z$.

On the other hand, $\Omega'$ is a closed form. To see this: Let
$i:M_t\rightarrow\mathcal{M}$ be the embedding of fibers, then
$i^*d\Omega'=d(i^*\Omega')=0$. It suffices to show $i_X d\Omega'=0$ since the restriction of $d\Omega'$ to the fiber is zero
already. Now we show $i_X d\Omega'=L_X\Omega'-di_X\Omega'=L_X\Omega'=0$. Notice that $\Omega'$ is determined by $\Theta|_{L_t}$ and the condition
$i_X\Omega'=0$. If we can show $\Theta|_{L_t}$ and $X$ are preserved by $X$-flow, then immediately we obtain $L_X\Omega'=0$ by uniqueness.
The fact $\Theta|_{L_t}$ is preserved follows $G^*\Theta$ is
constant along leaves and the fact $X$ is preserved follows $[X,\bar
X]=0$. So $\Omega'$ is closed form on $\mathcal{M}$, and
$i(h-\rho_{z\bar z}-\phi_{z\bar z})dzd\bar
z=\Omega'-\Omega-i\partial\bar\partial\phi$ is closed. This implies
$(h-\rho_{z\bar z}-\phi_{z\bar z})$ is just a function of $z$. Also,
since $\Omega'$ and $\Omega$ and $\phi$ are globally defined, so
$(h-\rho_{z\bar z}-\phi_{z\bar z})dzd\bar z$ is defined globally and
doesn't depend on the local representation. Therefore, the function
$h-\rho_{z\bar z}-\phi_{z\bar z}$ is globally defined, since
$dzd\bar z$ is defined on the whole disc. (Notice that the $z$
stands for a coordinate in a local product chart, so in different
product charts, $\phi_{z\bar z}$ is not the same though the function
$\phi$ is the same.)

Now let $H=h-\rho_{z\bar z}-\phi_{z\bar z}$. $H$ is defined
globally on $\pi\circ G(D\times U)$, but solely depends on $z\in D$.
One can solve the following equation on disc:
\begin{equation}
\partial_{z\bar z}\phi'=H
\end{equation}
with $\phi'=\phi_0-\phi$ on the $\partial D$. Now replace $\phi$ by
$\phi+\phi'$, then one get
$\Omega'=\Omega+i\partial\bar\partial\phi$ and $\phi=\phi_0$ on
$\partial D$. (Note that in different local charts, $(z,w)$ and
$(v,x)$ in $\mathcal{M}$, where $z,v$ project down to the same disc
variable. $\partial_{z\bar z}\phi'=\partial_{v\bar v}\phi'$ since
$\phi'$ is constant fiber-wise.) This finishes the proof of finding
potential $\phi$ if $U$ is sufficiently small.

Now for arbitrary $U$, one can always partition it into small open
balls $U_i$ which admit potential $\phi_i$. Let $\rho$ be a local
potential for the $\Omega$ on $\mathcal{M}$, then on the leaf

\begin{lem}$\Delta (\rho+\phi_i)=X\bar X(\rho+\phi_i)=0$
\end{lem}
\proof Let $f=\rho+\phi_i$,
\begin{eqnarray}
X\bar X f&=&X(\eta^{\bar\beta})f_{\bar\beta}+\partial\bar\partial f(X,\bar X)\\
&=&0
\end{eqnarray}
$\Box$

This implies $\Delta (\phi_i-\phi_j)=0$ on the leaf. Now with the
extra condition $\phi_i=\phi_j=\phi_0$ on the $\partial D$, it
implies $\phi_i=\phi_j$ on the intersection. The global potential is
immediately obtained from this. $\Box$

\begin{rem}The above correspondence is constructed only on simple
test configurations. In these configurations, central
fiber are smooth. However, we believe the techniques should work for some mild
singularities in the central fiber. 

Another point is that the correspondence has nothing to do with the
$C^*$ action.
\end{rem}

\section{Openness of super regular solution}
Using the correspondence in previous section, we can study
regularity of the solution $\phi$ by the associated holomorphic disc
family in $\mathcal{W}\rightarrow\mathcal{M}$.\footnote{However, the
existence so far only requires smoothness of total space.}
Donaldson's definition~\cite{DonaldsonSetup} of super regular discs
and the linearized model could be extended to our case as well. In
detail,
\begin{defn}
In the moduli map $G:D\times U\rightarrow\mathcal{W}$, a
disc $G(D,x)$ is called super regular at $z\in D$ if $d(\pi\circ
G_z)_x:TU\rightarrow TM$ is isomorphism. A disc $G(D,x)$ is called
super regular if it is super regular at every $z\in D$.
\end{defn}

\begin{defn}\label{SuperReg}A geodesic ray induced from a simple test configuration is called super regular if the disc family in $\mathcal{W}$ is super regular.  \footnote{i.e.: the solution is smooth regular to the Monge Ampere equation on the test configuration $\mathcal{M}$}.
\end{defn}

For a disc $G_x=G(\cdot,x)$ in the moduli map $G:D\times
U\rightarrow\mathcal{W}$, one can consider the holomorphic perturbation of
$G_x$ that satisfies the totally real boundary condition (the
boundary is in the $\Lambda_\phi$, i.e., the lifting of $M_t,
t\in\partial D$ by $\Omega+i\partial\bar\partial\phi$). Also, we normalize the perturbation such that it
preserves the projection property. In another word, $\pi\circ G:D\times U\rightarrow
D$ is identity on $D$ variable. The linearized problem is

\begin{thm} In the moduli map $G:D\times U\rightarrow\mathcal{W}$
corresponding to a smooth solution $\phi$, the linearized
perturbation equation for a disc $G(\cdot,x)$ is
\begin{eqnarray}
v&=&Su+A\bar u\texttt{   on   }\partial D \\
\bar\partial u&=&0\\
\bar\partial v&=&0
\end{eqnarray}
where $S$ and $A$ are maps from $\partial D$ to complex symmetric
matrices and positive hermitian matrices respectively. $u,v$ are
$C^n$ valued functions on $D$.
\end{thm}

\proof The idea is the same to Donaldson~\cite{DonaldsonSetup}. Trivialize the exact sequence $0\rightarrow (\pi\circ
G_x)^*(T^*\mathcal{M})\rightarrow G_x^*(T\mathcal{W})\rightarrow
(\pi\circ G_x)^*(T\mathcal{M})\rightarrow 0$. $\Box$

In~\cite{DonaldsonSetup}, it is showed that the problem is Fredholm and
the index is $2n$. Consequently, if the disc is regular in Fredholm sense, then $G:D\times U\rightarrow\mathcal{W}$ is indeed an open set
in the universal moduli space.

Regarding on the criterion of regularity of a disc, a modification of
Donaldson's argument leads to the following:

\begin{thm}If a disc is super regular at any point $p\in\partial D$,
then the disc is regular.
\end{thm}

\proof We look at the linearized model since the general case can be
simplified to the model.

First, define $\Omega(s_1,s_2)=u_1^tv_2-u_2^tv_1$. This is a
symplectic form for $s=(u,v)\in C^{2n}$. In particular, for
$s_1,s_2\in \ker\bar\partial_{S,A}$, $i\Omega(s_1(\tau),s_2(\tau))$
is real and independent of $\tau$. To see this, just notice that
$i\Omega(s_1,s_2)$ is holomorphic function and on $\partial D$,
$i\Omega(s_1,s_2)=i[u_1^t(Su_2+A\bar u_2)-u_2^t(Su_1+A\bar
u_1)]=i(u_1^tA\bar u_2-u_2^tA\bar u_1)$ is real.

The super regularity at $p\in\partial D$ means there are $2n$
elements $s_j=(u_j,v_j)\in \ker\bar\partial_{S,A}$ such that
$u_j(p)$ form a $R$-basis for $C^n$. By continuity, it implies
$u_j(\tau)$ form a $R$-basis for $C^n$ in a neighborhood $\tau\in
U_p$.

We claim $s_i(\tau)$ are generically $C$-linearly independent. It is
equivalent to claim $\det[s_j]_{1\leq j\leq 2n}$ has discrete zeros.
Notice $\det$ is holomorphic, so the zeros are either discrete or
the whole disc. Suppose it is the whole disc for contradiction. In
the neighborhood $U_p$, assume the maximal rank of $[s_j]_{1\leq
j\leq 2n}$ for $\tau\in U_p$ is achieved at $p$ without loss of
generality, and the rank is $k<2n$. Assume $s_1,s_2,...,s_k$ form a
basis for $span\{s_i\}$ at $p$, then near $p$,
$s_{k+1}=\sum\lambda_i s_i, 1\leq i\leq k$. $\lambda_i$ is
holomorphic, since it satisfies $\sum\lambda_is_i^ts_j=s_{k+1}^ts_j,
1\leq i,j\leq k$. In another word, it is obtained by solving the
holomorphic matrix equation $\lambda[s_i^ts_j]=s_{k+1}^ts_j$. Now
one finds holomorphic functions
$\lambda_1,...,\lambda_k,\lambda_{k+1}=-1,\lambda_{k+2}=0,...,\lambda_{2n}=0$
near $p$, such that $\sum \lambda_i s_i=0$. On the boundary
$\partial D$ near $p$,
\begin{equation}
0=\sum \lambda_jv_j=S(\sum\lambda_j u_j)+A(\sum\lambda_j\bar
u_j)=A(\sum\lambda_j\bar u_j)
\end{equation}
So $\sum\lambda_j\bar u_j=0$ and we also have $\sum\lambda_j u_j=0$,
so
\begin{equation}
\sum Im(\lambda_j)u_j=0=\sum Re(\lambda_j)u_j
\end{equation}
Since $u_j$ form $R$-basis near $p$, one has $\lambda_j=0$ on
$\partial D$ near $p$, which contradicts the choice of $\lambda_j$.
Therefore, the $\det[s_j]_{1\leq j\leq 2n}$ has discrete zero.

Now suppose the $\ker\bar\partial_{S,A}$ has dimension strictly
greater than $2n$. Then one can choose $s_0$ not in
$span\{s_i\},1\leq i\leq 2n$. Now in the $2n+1$ dimensional vector
space $span\{s_i\}$, $i\Omega$ as a skew form, must be singular. So
there is a vector $s\in span\{s_0,...,s_{2n}\}$ such that
$i\Omega(s,span\{s_1,...,s_{2n}\})=0$. Notice we proved
$s_1,...,s_{2n}$ form a $C$-basis generically, this implies $s=0$
generically on $D$. Thus it implies $s=0$, contradiction. $\Box$

In particular, since the holomorphic discs associated to smooth
solution $\phi$ are automatically super regular, above theorem
proves that they are all regular and the moduli space $M$ in
the map $G:D\times M\rightarrow\mathcal{W}$ is a compact connected
component of the universal moduli space. It readily implies the
following theorem.

\begin{thm}\label{Open}Openness:
If the equation
$(\Omega+i\partial\bar\partial\phi)^{n+1}=0,\phi=\phi_0$ on
$\partial\mathcal{M}$ admits a smooth solution $\phi$ with
$\Omega+i\partial\bar\partial\phi>0$ on fibers, then for any small
perturbation $\delta\phi_0\in C^\infty(\partial\mathcal{M})$, the
new boundary value problem still has smooth solution $\phi'$ which
is close to $\phi$ in $C^\infty(\mathcal{M})$ and
$(\Omega+i\partial\bar\partial\phi')>0$ on fibers.
\end{thm}

\proof We refer the proof to~\cite{DonaldsonSetup}, which
essentially asserts that compact families of regular normalized
discs are stable under small perturbations. $\Box$

\section{Geodesic ray from Toric degenerations}
\subsection{Basics of Toric degeneration}
For completeness, we describe Donaldson's construction of Toric degenerations~\cite{DonaldsonToric} in the following:

Let $P\subset R^n$ be a polytope associated to a toric variety $M$, for simplicity, assume $P$ is Delzant. Given a rational piece wise linear function $f$ on $M$, one associates with a polytope $\hat P=\{(x,y): x\in P, 0\leq y\leq K-f\}\subset R^{n+1}$, $K>\max f$. For simplicity, we assume $\hat P$ is Delzant and integral.

It is classical fact that $\hat P$ as above induces a toric variety $\mathcal{M}$ with a positive line bundle $\mathcal{L}$. Each integral point $p$ in $\hat P$ corresponds to a section $s_p$ of $\mathcal{L}\rightarrow\mathcal{M}$. The correspondence is compatible with addition  of integral points and multiplication of sections. In another word, if $p_1+p_2=p_3+p_4$, then $s_{p_1}s_{p_2}=s_{p_3}s_{p_4}$.   

One can view $\mathcal{M}$ as a sub-variety in $P^N$ by Kodaira embedding: $x\in\mathcal{M}, x\rightarrow [s_1(x):s_2(x):...:s_i(x)...]$, where $i$ runs through the integral points of $\hat P$. So $\mathcal{M}\subset P^N$ is defined by homogeneous equations $F(X_i)=0$. These equations are induced by the relations of $s_i$, or equivalently, by the relations of the integral points in $\hat P$.

There is a map $\pi: \mathcal{M}\rightarrow P^1$, defined by $\pi: x\rightarrow [s_p(x):s_q(x)]$ where $p=(t_1,...,t_n,t_{n+1}), q=(t_1,...,t_n,t_{n+1}+1)\in \hat P$. Also, there is a natural $C^*$ action on $\mathcal{M}$ from the torus $T^{n+1}=T^n\times C^*$. It transforms section $s_p$ to $t^ks_p$ where $p=(t_1,...,t_n,k)$. So the $C^*$ action can be lifted to $\pi:\mathcal{M}\rightarrow P^1$ by defining $t \circ [x:y]=[x:ty]$ on $P^1$.

The toric degeneration is just $\mathcal{M}-\pi^{-1}([1:0])$. The following example shows the construction in detail.

\ex:  Let $P=[0,2]\in R$ be the base polytope. $f=\max\{0,x-1\}$ is the piece wise linear function on $P$. $\hat P=([0,1]\times[0,1])\bigcup\{1\leq x\leq 2, x+y\leq 2\}$. 

Denote the integral points $X=(0,0), Y=(1,0), Z=(2,0), U=(0,1), V=(1,1)$. Then the toric degenerations is the sub-variety in $P^4$ defined by
\begin{equation}
XZ=Y^2, XV=UY
\end{equation}
The $C^*$ action on $\mathcal{M}$ is $t:[X:Y:Z:U:V]\rightarrow[X:Y:Z:tU:tV]$. Notice that in order to get nontrivial test configuration, we only consider the part $\mathcal{M}-\pi^{-1}([1:0])$. In another word, we consider the asymptotic direction when $t\rightarrow \infty$ on $C^*$.

The central fiber is defined by $[Y:V]=[0:1]$. It is the toric variety associated to the segment $y=1,x\in[0,1]$ and $x\in[1,2],x+y=2$. Geometrically, the central fiber is the union of two $P^1$ which intersect at one point.

Notice that the ambient space $\mathcal{M}$ is smooth here, so the induced geodesic ray has ambient bounded geometry automatically.

\subsection{Explicit calculation of the $C^{1,1}$ geodesic ray}
We calculate the induced geodesic ray of previous example. The idea is to first calculate the geodesic segment connecting the fiber at $[1:1]$ to the fiber at $[1:e^t], t\in R\times S^1$, and then take the limit of these segments when $t\rightarrow \infty$.

Equipped with the natural background metric of $P^4$, the fiber at $w=[1:e^t]\in P^1$ has metric potential $\frac{1}{2}\log(|X|^2+|Y|^2+|Z|^2+|U|^2+|V|^2)$. Pull this metric to the fixed fiber $M$ at $w=[1:1]\in P^1$, the potential becomes
\begin{equation}
\frac{1}{2}\log(|X|^2+|Y|^2+|Z|^2+e^{2t}|U|^2+e^{2t}|V|^2)
\end{equation}
Since the fiber $M$ is at $[1:1]$, so $Y=V,X=U$. After proper normalization, the potential is
\begin{equation}
\frac{1}{2}\log(|X|^2+|Y|^2+(e^{2t}+1)^{-1}|Z|^2)
\end{equation}

Now we calculate the geodesic segment connecting these two metrics. 

Choose $[A,B]$ as standard $P^1$ coordinate on $M$, so $X=B^2, Y=AB, Z=A^2$. Using $C^*=R\times S^1$ coordinate of $P^1$, $A=e^y,B=1, y\in R\times S^1$, and the metric potential is 
\begin{equation}
h_{0,t}=\frac{1}{2}\log(1+e^{2y}+e^{4y}(e^{2t}+1)^{-1})
\end{equation}
One can verify the legendre transform of $h_{0,t}$ maps $R$ to $(0,2)$ for each fixed $t$. 

Notice that in polytope representation, the geodesic is just a straight line of convex functions. 
Now by straightforward calculation, one just computes the two end points associated to the two metrics in polytope representation and then take the linear interpolation. Passing to limit, one gets
the $C^{1,1}$ ray in polytope representation
\begin{equation}
u_t=u_0+t\max(0,x-1)  \texttt{       ,      }   t\in [0,\infty)
\end{equation}

In the standard picture of $M\times [0,\infty)$, we transform the $u_t$ by Legedre transform and get the potential
\begin{equation}
h_t(y)=\left\{
\begin{array}{ll}
  h_0(y)  & y<\frac{\log 2}{4}   \\
  h_0(\frac{\log 2}{4})+y-\frac{\log 2}{4}&    \frac{\log 2}{4}<y<\frac{\log 2}{4}+t  \\
  h_0(y-t)+t& \frac{\log 2}{4}+t<y   
\end{array}
\right.
\end{equation}

One can verify that $h_t-h_{0,t}$ is uniformly bounded. This confirms that the geometric ray is parallel to the algebraic ray.\\

It is natural to extend above result to general toric degenerations. 
\begin{thm}Let $\mathcal{M}$ be a toric degeneration with extremal piece wise linear function $f$. Suppose the ambient polytope $\hat P$ is delzant. Then the induced geodesic ray is $u=u_0+tf$ in polytope representation.
\end{thm}
\proof We only give a sketch of the proof. Because of the uniqueness of geodesic ray in a fix direction, it suffices to show the ray $u=u_0+tf$ is parallel to the algebraic ray. In another word, it suffices to show $h_t-h_{0,t}$ is uniformly bounded. $h_t$ is the geodesic ray in standard product presentation. $h_{0,t}$ is the algebraic potential.

If we take the fiber $M_1$ at $[1,1]$ for standard model, then $h_{0,t}=\frac{1}{2}\log\Sigma e^{2k_i}|X_i|^2$, with $X_i$ are integral points in the base polytope. Using $C^*=R\times S^1$ coordinates, one can write $X_i=\exp(\Sigma d_jy_j)$, and $X_i=(d_1,d_2,...,d_n)$. 

The calculation of $h_t$ is the same as in the example. It is straightforward but lengthy to check the difference $h_t-h_{0,t}$ is uniformly bounded. $\Box$.\\

These geodesic rays show some bad regularity. In general, they behave like the following: First, they break the manifold $M$ into several pieces. As time evolves, they will tear these pieces apart, but keep metrics on each part. The space between the teared parts has vanished metric. In particular, one can verify that the 2nd derivative of these rays are piece wise smooth function on fibers. At the broken points, these 2nd derivatives have jumps, so there is no global $C^3$ bound for the geodesic ray potential. \\

For toric varieties, there has been extensive literature in extremal metrics. Abreu\cite{Abreu} initiated to study complex geometry on toric variety by symplectic coordinates. Afterwards, there has been much work in extremal metrics on toric variety, c.f. Donaldson\cite{DonaldsonToric}, Zhou-Zhu\cite{Zhu},Gabor\cite{Gabor}.

\end{document}